\documentclass[10pt]{article}
\usepackage{amssymb, amsmath, url, graphicx, setspace, geometry}
\usepackage{calrsfs}
\usepackage{wasysym, xcolor}

\title{\bf Minimum $\mathcal{F}$-covers: the supersolvable and metabelian cases}
\author{Mihai-Silviu Lazorec}

\begin{document}

\maketitle

\begin{abstract}
Given a set $\mathcal{F}$ of finite groups, it is said that a group $G$ is an $\mathcal{F}$-cover if every group in $\mathcal{F}$ is isomorphic to a subgroup of $G$. Moreover, $G$ is a minimum $\mathcal{F}$-cover if there is no $\mathcal{F}$-cover whose order is less than $|G|$. In [Cameron P. J., et al., Minimal cover groups, J. Algebra 660 (2024)], the authors pose the following question: For which classes $\mathcal{X}$ of groups, closed under taking subgroups and direct products, is it true that, if $\mathcal{F}$ is a set of $\mathcal{X}$-groups, then there is a minimum $\mathcal{F}$-cover which is an $\mathcal{X}$-group? In this paper, we give a negative answer in two cases: $\mathcal{X}\in \{``supersolvable", ``metabelian"\}.$  
\end{abstract}

\noindent{\bf MSC (2020):} Primary 20F16; Secondary 20D40, 20E99.

\noindent{\bf Key words:} cover of a set of groups, supersolvable group, metabelian group.  

\section{Introduction}

All groups considered in this paper are finite. $\mathcal{F}$ always denotes a finite set of  groups. The order of an element $x$ of a group, the cyclic group of order $n$ ($n\geq 2$), the dihedral group of order $2n$ ($n\geq 3$) and the Schur multiplier of a group $G$ are denoted by $o(x)$, $C_n$, $D_{2n}$ and $M(G)$, respectively. 

We say that a group $G$ is an $\mathcal{F}$-cover if every group in $\mathcal{F}$ is isomorphic to a subgroup of $G$. Inspired by Cayley's theorem stating that any group of order $n$ is isomorphic to a subgroup of $S_n$, the $\mathcal{F}$-cover concept is  introduced and studied by Cameron, Craven, Dorbidi, Harper and Sambale in \cite{1}. Some of the problems that are approached there and the  corresponding results mainly refer to two types of $\mathcal{F}$-covers: minimal and minimum. Assuming that $G$ is an $\mathcal{F}$-cover, it is said that this cover is:
\begin{itemize}
\item[--] minimal if no proper subgroup of $G$ is an $\mathcal{F}$-cover;
\item[--] minimum if there is no $\mathcal{F}$-cover whose order is less than $|G|$. 
\end{itemize}
It is clear that any minimum $\mathcal{F}$-cover is also minimal. 

One of the main results in \cite{1} is stated as follows:\\

\textbf{Theorem 1.} \textit{Let $\mathcal{F}=\{ C_q, C_r\}$, where $q, r$ are  prime numbers such that $q<r$. Then $\mathcal{F}$ has finitely many minimal covers if and only if $q=2$ and $r$ is a Fermat prime. More exactly, these covers are: $C_{2r}, D_{2r}$ and $C_2^{2a}\rtimes C_r$, where $r=2^a+1$.}\\
 
Hence, if $\mathcal{F}$ is formed of two cyclic groups of prime order, then Theorem 1 provides an answer to the following question: \textit{For which sets $\mathcal{F}$ are there only finitely many minimal $\mathcal{F}$-covers up to isomorphism?} (see Question A in \cite{1}).

A further research direction that is outlined in \cite{1} is Question 8.2: \textit{For which classes $\mathcal{X}$ of groups, closed under taking subgroups and direct products, is it true that, if $\mathcal{F}$ is a set of $\mathcal{X}$-groups, then there is a minimum $\mathcal{F}$-cover which is an $\mathcal{X}$-group?} Regarding this direction, the authors prove the following result (see Theorem 7.3 in \cite{1}) that provides an affirmative answer to Question 8.2 when $\mathcal{X}=``nilpotent"$.\\
 
\textbf{Theorem 2.} \textit{Let $\mathcal{F}$ be formed of nilpotent groups. Then there is a minimum $\mathcal{F}$-cover which is nilpotent.}\\

A way to construct a minimum $\mathcal{F}$-cover in the nilpotent case is shared in the proof of Theorem 7.3 of \cite{1}. Example 7.8 of the same paper shows that if $\mathcal{F}=\{ D_{10}, A_4\}$, then $A_5$ is the unique minimum $\mathcal{F}$-cover. This shows that the answer to Question 8.2 is negative if $\mathcal{X}=``solvable"$. If $\mathcal{X}=``abelian"$, Question 8.2 can be reduced to the case of abelian $p$-groups. At the moment of writing this paper, it is unknown if starting with $\mathcal{F}$ formed of abelian $p$-groups, then there is a minimum $\mathcal{F}$-cover which is an abelian $p$-group. Still, if this is true, a way to construct this cover is outlined in the paragraph that precedes Proposition 7.4 of \cite{1}. 

Even though by taking $\mathcal{X}=``cyclic"$, we get a class of groups that is closed only under taking subgroups, we think that it is worth mentioning that if $\mathcal{F}$ contains only cyclic subgroups, we can always find a minimum $\mathcal{F}$-cover which is cyclic. Theorem 7.1 of \cite{1} outlines this cover. This result is stated below.\\

\textbf{Theorem 3.} \textit{Let $\mathcal{F}=\{ C_{n_1}, C_{n_2}, \ldots, C_{n_k}\}$ and $N=lcm(n_1, n_2, \ldots, n_k)$. Then $C_N$ is a minimum $\mathcal{F}$-cover.}\\  

In this paper, we provide a negative answer to Question 8.2 in two more specific cases: $$\mathcal{X}\in \{ ``supersolvable", ``metabelian"\}.$$ We outline infinitely many sets $\mathcal{F}$ formed of $\mathcal{X}$-groups such that they have a unique minimum $\mathcal{F}$-cover which is not an $\mathcal{X}$-group. More exactly, we prove the following two results:\\

\textbf{Theorem 4.} \textit{Let $p\geq 5$ be a prime number and let $\mathcal{F}_p=\{  C_3\rtimes S_3, C_{4p}\}$ be a set of supersolvable groups. Then the unique minimum $\mathcal{F}_p$-cover is the non-supersolvable group $C_p\times (C_3^2\rtimes C_4)$.}\\

\textbf{Theorem 5.} \textit{Let $p$ be an odd prime such that $p\neq 5$ and let $\mathcal{F}_p=\{ D_{10}, C_p\times A_4\}$ be a set of metabelian groups. Then the unique minimum $\mathcal{F}_p$-cover is the non-metabelian group $C_p\times A_5$.}\\

Regarding the statement of Theorem 4, we note that $C_3\rtimes S_3$ is SmallGroup(18,4) and $C_3^2\rtimes C_4$ is SmallGroup(36, 9) in GAP's \cite{4} small groups library. The proofs of the two results above are given in the following section.

We end the introduction by pointing out that if ``minimum" is replaced with ``minimal" in Question 8.2, its corresponding answer is always affirmative by Proposition 2.6 of \cite{1}.      
    
\section{Proofs of main results}

As preliminary results, we recall Proposition 2.1 (\textit{a}) of \cite{1} and Theorem 2.1 (\textit{a}) of \cite{3}. The first one is a   consequence of Lagrange's theorem and it gives us a clue on the order of a minimum $\mathcal{F}$-cover. The second one provides information on the structure of a group having a cyclic normal subgroup with specific properties.\\

\textbf{Lemma 2.1.} \textit{Let $G$ be an $\mathcal{F}$-cover. Then lcm$\{|F| \mid F\in\mathcal{F} \}$ divides $|G|$.}\\

\textbf{Lemma 2.2.} \textit{Let $G$ be a group, $H$ be a cyclic normal subgroup of $G$ and $K$ be a simple non-abelian group. If $\frac{G}{H}\cong K$ and gcd$(|H|, |M(K)|)=1$, then $G\cong H\times K$.}\\

\textit{Proof of Theorem 4.} Let $G$ be a minimum $\mathcal{F}_{p}$-cover. By Lemma 2.1, we have that $36p\mid |G|$, establishing $36p$ as the minimum possible order for any $\mathcal{F}_{p}$-cover. In light of the existence of the $\mathcal{F}_{p}$-cover $C_{p}\times(C_{3}^{2}\rtimes C_{4})$ of order $36p$ given in the statement of the theorem, we may assume that $|G|=36p$. It therefore suffices to show that this cover is unique. 

If $p\in \{5, 7, 11, 13, 17\}$, then one can use GAP to check that $G\cong C_p\times (C_3^2\rtimes C_4)$ is the unique minimum $\mathcal{F}_p$-cover. 

In what follows, we assume that $p\geq 19$. By using Sylow's theorems, it is easy to check that $G$ has a unique Sylow $p$-subgroup $H\cong C_p$. Hence, $H$ is normal in $G$ and since its order is coprime with the order of the quotient $K=\frac{G}{H}$, it follows that $$G\cong H\rtimes K$$ by the Schur-Zassenhaus theorem. Let $\pi:G\longrightarrow K$ be the surjective projection homomorphism. Since $G$ is an $\mathcal{F}_p$-cover, we know that $G$ has subgroups isomorphic to $F_1=C_3\rtimes S_3$ (SmallGroup(18, 4)) and $F_2=C_{4p}$, respectively. Consequently, $G$ has a subgroup $L_1\cong S_3$ and a subgroup $L_2\cong C_4$. Let $\pi_{/L_1}:L_1\longrightarrow K$ be the restriction of $\pi$ to $L_1$. We easily obtain that its kernel is
$$ker(\pi_{/L_1})=H\cap L_1=\{ 1\}.$$
Hence, $\pi_{/L_1}$ is actually a group isomorphism from $L_1$ to $\pi_{/L_1}(L_1)$, so $K$  has a subgroup isomorphic to $S_3$. By following a similar reasoning, one deduces that $K$ also has a  subgroup isomorphic to $C_4$. For the remainder of the proof, we may choose $L_1$ and $L_2$ to be these specific subgroups of $K$. By inspecting the list of groups of order 36, the only possible choice is $K\cong C_3^2\rtimes C_4$ (SmallGroup(36, 9)). Hence,
$$G\cong H\rtimes K\cong C_p\rtimes (C_3^2\rtimes C_4).$$
Let $f:K\longrightarrow Aut(H)$ be the group homomorphism that induces the above semidirect product. We are going to show that $f$ is trivial. It is known that $ker(f)=C_{K}(H)$ (see, for instance, Exercise 1, Section 5.5 of \cite{2}). Without loss of generality, we may assume that $F_2= L_2\times H$. Then $L_2\subseteq C_K(H)$. Let $N\cong C_3^2$ be the Sylow 3-subgroup of $K$. Since $C_K(H)$ and $N$ are normal subgroups of $K$, we have $M=C_K(H)\cap N\triangleleft K$. The group $K$ does not have normal subgroups of order 3, so we deduce that $M=\{1\}$ or $M=N$. In the first case, $C_K(H)N$ would be a subgroup of order $9|C_K(H)|$ of $K$. By taking into account that $L_2\subseteq C_K(H)$, it follows that $C_K(H)=L_2\cong C_4$. It follows that $K\cong C_3^2\times C_4$, a contradiction. Hence, $M=N$ and this implies that $N\subseteq C_K(H)$. Since $N\cap L_2=\{ 1\}$ and $NL_2\subseteq C_K(H)$, we get that $C_K(H)$ is a subgroup of $K$ with $|C_K(H)|\geq 36$. It follows that $C_K(H)=K$, so $f$ is trivial. Then,
$$G\cong H\times K\cong C_p\times (C_3^2\rtimes C_4),$$
is the unique minimum $\mathcal{F}_p$-cover.

For any prime $p\geq 5$, $G$ is an $\mathcal{F}_p$-cover which is a solvable group, but it is not supersolvable since its subgroup $K$ is a non-supersolvable group. The proof is complete.
\hfill\rule{1,5mm}{1,5mm}\\

With the notations established in the statement of Theorem 4, we mention that if $p=2$, then the unique minimum $\mathcal{F}_2$-cover is a non-supersolvable Frobenius group of order 72: $F_9=C_3^2\rtimes C_8$ (SmallGroup(72, 39)). If $p=3$, there are 4 minimum $\mathcal{F}_3$-covers, all of them being supersolvable groups: 
\begin{itemize}
\item[--] $C_3^2\rtimes (C_2\times C_4)$ (SmallGroup(72, 21));
\item[--] $C_3\rtimes D_{24}$ (SmallGroup(72, 23));
\item[--] $C_4\times (C_3\rtimes S_3)$ (SmallGroup(72, 32));
\item[--] $C_{12}\rtimes S_3$ (SmallGroup(72, 33)).
\end{itemize}

If $\mathcal{F}$ is formed of supersolvable groups, it may happen that there are at least 2 minimum $\mathcal{F}$-covers, all of them being non-supersolvable. For instance, if $\mathcal{F}=\{ C_4\rtimes C_4, C_3\rtimes S_3\}$, there are 3 minimum $\mathcal{F}$-covers and all are non-supersolvable: SmallGroup(144, 116), SmallGroup(144, 120), SmallGroup(144, 133). These 3 groups are some semidirect products of $C_3^2$ and $C_4\rtimes C_4$. The same 3 groups are also the unique minimum $\mathcal{F}$-covers of $\mathcal{F}=\{ C_4\rtimes C_4, C_2\times (C_3\rtimes S_3)\}$. This shows that two different sets of supersolvable groups may have the same minimum covers, all of them being non-supersolvable. 

A second example of infinitely many sets $\mathcal{F}$ formed of supersolvable groups such that they have a unique minimum $\mathcal{F}$-cover which is not a supersolvable group can be given. More exactly, if $p\geq 5$ is a prime and $\mathcal{F}_p=\{ C_3\rtimes S_3, C_{8p}\}$ is a set of supersolvable groups, its unique minimum $\mathcal{F}_p$-cover is the non-supersolvable group $C_p \times F_9$.  An argument would follow similar ideas that are part of the proof of Theorem 4.\\

\textit{Proof of Theorem 5.} Let $G$ be a minimum $\mathcal{F}_p$-cover. By a similar argument as the one given in the proof of Theorem 4, we may assume that $|G|=60p$ and it suffices to prove uniqueness.

If $p\leq 59$, then one can check that the theorem holds via GAP. If $p\geq 61$, then $G$ has a cyclic normal subgroup $H\cong C_p$. Let $K=\frac{G}{H}$. By following a part of the reasoning that was done in the proof of Theorem 4, we deduce that $K$ has  subgroups isomorphic to $D_{10}$ and $A_4$, respectively. Up to isomorphism, $A_5$ is the unique group of order 60 satisfying such properties, so $K\cong A_5$. Then $K$ is a simple non-abelian group and, since gcd$(|H|, |M(K)|)=$ gcd$(p, 2)=1$, we obtain that $$G\cong H\times K\cong C_p\times A_5,$$
by Lemma 2.2. 

Obviously, $G$ is a non-solvable group, so it is also non-metabelian. The proof is complete. 
\hfill\rule{1,5mm}{1,5mm}\\  

By following the notations established in the statement of Theorem 5, we also study the cases $p\in\{2, 5\}.$ If $p=2$, then there are  two minimum $\mathcal{F}_2$-covers: $C_2\times A_5$ and $D_{10}\times A_4$. If $p=5$, there are also two minimum $\mathcal{F}_5$-covers: $C_5\rtimes S_4$ (SmallGroup(120, 38)) and $D_{10}\times A_4$. For both values of $p$, the first corresponding minimum $\mathcal{F}_p$-cover is a non-metabelian group, while the second is metabelian.

As in the supersolvable case, we mention that for specific choices of $\mathcal{F}$ formed of metabelian groups, we can get at least 2 minimum $\mathcal{F}$-covers, all of them being non-metabelian. As an example, if one  chooses $\mathcal{F}=\{D_8, C_3\times A_4\}$, then the only minimum $\mathcal{F}$-covers are $C_3\times S_4$ and $C_3\rtimes S_4$ (SmallGroup(72, 43)). Also, a non-metabelian group can be the unique minimum $\mathcal{F}$-cover for different choices of $\mathcal{F}$ formed of metabelian groups. For instance, $S_3^2\rtimes C_2$ (SmallGroup(72, 40)) is the unique minimum $\mathcal{F}$-cover if $\mathcal{F}\in\{ \{D_8, C_3^2\rtimes C_4\}, \{D_8, S_3^2 \}\}$.

Finally, we point out infinitely many sets $\mathcal{F}$ formed of metabelian groups such that they have only two minimum $\mathcal{F}$-covers, both of them being non-metabelian. More precisely, if $p\geq 5$ is a prime and $\mathcal{F}_p=\{ D_8, C_p\times A_4\}$, its only minimum $\mathcal{F}_p$-covers are $C_p\times S_4$ and $C_p\rtimes S_4$. If $p\leq 23$, this can be checked via GAP. If $p\geq 29$, as in the proofs of Theorems 4 and 5, we get that a minimum $\mathcal{F}_p$-cover $G$ is isomorphic to $H\rtimes K$, where $H\cong C_p$ and $K$ is a group of order 24 having subgroups isomorphic to $D_8$ and $A_4$, respectively. It follows that $K\cong S_4$, so $G\cong C_p\rtimes S_4$. If $f: S_4\longrightarrow Aut(C_p)$ is the group homomorphism that induces the previous semidirect product, we deduce that
$$\frac{S_4}{ker(f)}\cong im(f)\leq Aut(C_p)\cong C_{p-1}.$$
Hence, $\frac{S_4}{ker(f)}$ is an abelian group and this leads us to $S_4'=A_4\leq ker(f)$. Therefore, $f$ induces a unique nontrivial semidirect product or a direct product. Consequently, $G\cong C_p\rtimes S_4$ or $G\cong C_p\times S_4$. 

\bigskip\noindent {\bf Acknowledgements.} The author is grateful to the reviewer for his/her remarks which improved the initial version of the paper. 
  
\bigskip\noindent {\bf Declarations}

\bigskip\noindent {\bf  Funding.} The author did not receive support from any organization for the submitted work.

\bigskip\noindent {\bf Conflicts of  interests.} The author declares that there is  no conflict of interest.

\bigskip\noindent {\bf  Data availability statement.} The manuscript has no associated data.

\vspace*{3ex}
\small
\hfill
\begin{minipage}[t]{7cm}
Mihai-Silviu Lazorec \\
Faculty of  Mathematics \\
"Al.I. Cuza" University \\
Ia\c si, Romania \\
e-mail: {\tt silviu.lazorec@uaic.ro}
\end{minipage}
\end{document}